\documentclass[10pt,a4paper]{amsart}
\usepackage[utf8]{inputenc}
\usepackage[T1]{fontenc}
\usepackage[english]{babel}
\usepackage{amsmath,hyperref}
\usepackage[all]{xy}

\usepackage[active]{srcltx}
\usepackage{enumerate}
\usepackage[parfill]{parskip}

\title[Homologically optimal categories of sequences]{Homologically optimal categories of sequences\\ lead to $N$-complexes}
\author{Djalal Mirmohades}
\address{Uppsala University, Sweden}
\email{math@djalal.com}

\newtheorem{theorem}{Theorem}
\newtheorem{proposition}[theorem]{Proposition}
\newtheorem{corollary}[theorem]{Corollary}
\newtheorem{lemma}[theorem]{Lemma}
\theoremstyle{definition}
\newtheorem{definition}[theorem]{Definition}

\newcommand{\Z}{\mathbb{Z}}
\newcommand{\I}{\mathrm{I}}
\newcommand{\F}{\mathrm{F}}
\newcommand{\G}{\mathrm{G}}
\newcommand{\T}{\mathrm{T}}
\newcommand{\Q}{\mathrm{Q}}
\newcommand{\R}{\mathrm{R}}
\newcommand{\E}{\mathrm{E}}
\newcommand{\HH}{\mathrm{H}}
\newcommand{\Id}{\mathrm{Id}}
\newcommand{\A}{\mathcal{A}}
\newcommand{\B}{\mathcal{B}}
\newcommand{\C}{\mathcal{C}}
\newcommand{\D}{\mathcal{D}}

\font\sc=rsfs10
\newcommand{\DD}{\sc\mbox{D}\hspace{1.0pt}}

\newcommand{\eps}{\varepsilon}
\newcommand{\dd}{[\ker d^N]}
\newcommand{\im}{\operatorname{im}}
\newcommand{\Seq}{\operatorname{Seq}}
\newcommand{\Mor}{\operatorname{Mor}}
\newcommand{\Com}[1]{\mathrm{Com}_{#1}}

\begin{document}

\begin{abstract}
We study the category of $\mathbb{Z}$-indexed sequences over an abelian category and certain generalized homology functors for this category of sequences which are indexed by positive integers $a$ and $b$. By looking at the corresponding derived  category, we show that there is an ``optimal'' subcategory of sequences for every choice of our generalized homology functors, namely, the category of $N$-complexes (sequences for which the differential $d$ satisfies $d^N = 0$) where $N = a + b$. In this optimal case we show that our homology functors reduce to Kapranov's homology functors $\operatorname{ker} d^a / \operatorname{im} d^b$.
\end{abstract}

\maketitle


\section{Introduction}\label{s1}

Homological algebra traditionally studies chain complexes over an abelian category $\A$, that is seqeunces of
objects in $\A$ equipped with a differential $d$ satisfying $d^2 = 0$. Associated to this setup is the classical
homology functor $\operatorname{ker} d / \operatorname{im} d$.

More generally one can study the so-called $N$-complexes for $N\in\{2,3,4,\dots\}$, that is seqeunces of
objects in $\A$ equipped with a differential $d$ satisfying $d^N = 0$. In \cite{kapr91} Kapranov generalizes
some of the tools of homological algebra for such $N$-complexes, inspired by their application in quantum groups.
Since then, various aspects of the homological algebra of $N$-complexes were studied by various authors, see
\cite{DV1,DV2,DV3,DV4,CSW,Ti,Gi,HK} and references therein. 

In this paper, for  a fixed generalized homology functor of the form
\begin{equation}\label{gen-homology}
\HH := \ker d^a \big/ \left( \ker d^a \cap \im d^b \right) = \left( \ker d^a + \im d^b \right) \big/ \im d^b
\end{equation}
we ask the  following question: Which ``nice'' subcategory $\B$ of the category of sequences is the ``minimal''
one with the property that  the corresponding derived category $\DD_{\HH}(\B)$ with respect to $\HH$ inherits
from $\A$ the property of having enough projectives? One of our results, combined in 
Theorem~\ref{proj} and Corollary~\ref{enough_proj}, 
asserts that the category $\Com{N}(\A)$ of $N$-complexes of $\A$ has these properties for $N=a+b$.

In fact, we even show more. Namely, in Theorem~\ref{theorem_eq} and the subsequent proposition we prove that the category $\Com{N}(\A)$ is ``optimal'' in the sense that for any category $\B$ containing $\Com{N}(\A)$ there is an equivalence
$$\DD_{\HH}(\B) \simeq \DD_{\HH}(\Com{N}(\A))$$
of the corresponding derived categories with respect to our choice of $\HH$.
In the category $\Com{N}(\A)$ we have $\im d^b \subset \ker d^a$ because $d^a d^b = d^N = 0$, consequently $\HH$ in \eqref{gen-homology} reduces to Kapranov's homology $\ker d^a / \im d^b$.
In particular, if we choose $a = b = 1$ we end up with $\Com{2}(\A)$, the category of ordinary chain complexes together with classical notions of homology and derived category.

The paper is organized as follows: In Section~\ref{s2} we collect basic definitions on $N$-complexes and
their homology. The statement about derived equivalence is proved in Section~\ref{s5}. The statement
on projective resolutions in proved in Section~\ref{s6}.

This paper is an adaptation of a part of author's Master Thesis \cite{Mi}.

\section{Preliminaries}\label{s2}

\subsection{\texorpdfstring{Sequences and $N$-complexes}{Sequences and N-complexes}}\label{s2.1}

Given an abelian category $\A$, define the category $\Seq(\A)$ as follows:

An \textbf{object} $C$ of $\Seq(\A)$, called \textit{a sequence}, is a collection of objects 
$\{C_i\}_{i \in \Z}$ together with a collection of morphisms $\{d_i : C_i \to C_{i+1} \}_{i \in \Z}$ in $\A$ 
(called the {\em differential} of $C$ and denoted by $d$) which is usually depicted as follows:
\[ \xymatrix@C=30pt{ 
\cdots \: \ar[r]^{d_{-2}}  & C_{-1} \ar[r]^{d_{-1}} & C_0 \ar[r]^{d_{0}} & C_1 \ar[r]^{d_{1}} & \: \cdots .
} \]

For $C,D\in \Seq(\A)$, a \textbf{morphism} $f\in \mathrm{Hom}_{\Seq(\A)}(C,D)$ is a commuting diagram over $\A$
of the following form:
\[ \xymatrix@C=30pt{ 
\cdots \: \ar[r]  & C_{-1} \ar[r] \ar[d]_{f_{-1}} & C_0 \ar[r] \ar[d]_{f_{0}} & C_1 \ar[r] \ar[d]_{f_{1}} & \: \cdots \phantom{.} \\
\cdots \: \ar[r]  & D_{-1} \ar[r] & D_0 \ar[r] & D_1 \ar[r] & \: \cdots .
} \]

The category $\Seq(\A)$ is the functor category from the path category of the graph
\[\xymatrix@C=36pt{ 
\cdots \: \ar[r] & \text{-}\mathtt{1} \ar[r] & \mathtt{0}\ar[r] & \mathtt{1} \ar[r] & \: \cdots \phantom{.}
}\]
to $\A$, in particular, it is abelian because $\A$ is abelian.

For $N\in\{1,2,3,\dots\}$ we define the category $\Com{N}(\A)$ of {\em $N$-complexes} over $\A$ as the full 
subcategory of $\Seq(\A)$ which consists of all object which satisfy the condition 
$d_{i+N-1}\cdots d_{i+1}d_i=0$ for all $i\in\mathbb{Z}$. As usual, we will just write $d^N=0$ for the latter
condition.

\subsection{\texorpdfstring{Homology for $N$-complexes}{Homology for N-complexes}}\label{homology}

For 2-complexes, which are also usually called \textit{chain complexes}, we have the inclusion 
\[\xymatrix{
\im d \: \ar@<-1pt>@{^{(}->}[r] & \ker d
}\]
and there is a functor $\Com{2}(\A) \to \Com{1}(\A)$, called \textit{homology}, given by
\begin{equation}\label{eq-homology_d_d}
\ker d \big/ \im d.
\end{equation}
For $N$-complexes, the differential induces the following poset of non-trivial inclusions:
\begin{equation}\label{eq-inclusion}
\xymatrix@=0pt@!R=70pt@!C=35pt{ 
& & & & \ker d^{N-1} & &  \\
& & & \ker d^{N-2} \ar@{^{(}->}[ur] & & \im d \ar@{_{(}->}[ul] & \\
& & \ker d^2  \ar@{.>}[ur] & & \im d^2 \ar@{_{(}->}[ul] \ar@{^{(}->}[ur] & & \\
& \ker d \ar@{^{(}->}[ur] & & \im d^{N-2} \ar@{_{(}->}[ul] \ar@{.>}[ur] & & & \\
 & & \im d^{N-1} \ar@{_{(}->}[ul] \ar@{^{(}->}[ur] & & & &
}
\end{equation}

Our goal is to study various generalizations of the cvlassical homology given by \eqref{eq-homology_d_d} 
to $N$-complexes. To define such homologies we can use ingredients shown on the diagram 
\eqref{eq-inclusion} and also their various sums and intersections. For example,
the principle of ``throwing away as much information as possible'' results in the homology 
\begin{equation}\label{eq-homology_d_im_d}
\ker d \big/ \left( \ker d \cap \im d \right).
\end{equation}
On the other hand, the principle of ``keeping as much information as possible'' results in the homology 
\begin{equation}\label{eq-homology_d_N-1_d_N-1}
\ker d^{N-1} \big/ \im d^{N-1}.
\end{equation}
In some sense, the latter does preserve more information than, say 
\begin{equation}\label{eq-homology_d_N-1_d}
\ker d^{N-1} \big/ \im d
\end{equation}
when $N > 2$. However, one should be aware of the fact that the homology \eqref{eq-homology_d_N-1_d} 
cannot be recovered only knowing the homology \eqref{eq-homology_d_N-1_d_N-1}. 
To see this, consider the following pair of $3$-complexes of abelian groups:
\[ \xymatrix@R=10pt{
\cdots \: \ar[r]  & \Z/8\Z \ar[r]^2 & \Z/8\Z \ar[r]^2 & \Z/8\Z \ar[r]^2 & \Z/8\Z \ar[r] & \: \cdots \\
\cdots \: \ar[r]  & 2\Z/4\Z \ar[r]^0 & 2\Z/4\Z \ar[r]^0 & 2\Z/4\Z \ar[r]^0 & 2\Z/4\Z \ar[r] & \: \cdots
} \]
It is easy to check that the homology $\ker d^{2} / \im d^{2}$ from  \eqref{eq-homology_d_N-1_d_N-1} has 
identical values for both of them, while the homology  $\ker d^{2} / \im d$ from \eqref{eq-homology_d_N-1_d} 
has different values on these two complexes.

Now fix $N \geq 1$. For positive integers $a, b$ such that $a + b \geq N$, define the homology 
functor $\HH^{(a, b)}_j : \Com{N}(\A) \to \A$, where $j \in \Z$, as follows:
\begin{equation}\label{eq-homology_d_a_d_b}
\HH^{(a, b)}_j := \ker d^a \big/ \im d^b \: \big(\mbox{evaluated at position $j$}\big).
\end{equation}

\begin{proposition}\label{prop11}
Equation \eqref{eq-homology_d_a_d_b} indeed defines a functor. 
\end{proposition}

\begin{proof}
Consider the functor
$$\mathrm{S}^{(a, b)} : \Com{N}(\A) \to \Com{2}(\A)$$
which maps an $N$-complex 
\begin{equation}\label{eq2.1}
\xymatrix@C=30pt{
\: \phantom{\mathrm{S}^{(a, b)}} C: & \cdots \: \ar[r]^d  & C_{-1} \ar[r]^d & C_0 \ar[r]^d & C_1 \ar[r]^d & \: \cdots \phantom{.}
} 
\end{equation}
to the following $2$-complex:
\[\xymatrix@C=30pt{ 
\: \mathrm{S}^{(a, b)} C: & \cdots \: \ar[r]^{d^a} & C_{-b} \ar[r]^{d^b} & C_0 \ar[r]^{d^a} & C_{a} \ar[r]^{d^b} & \: \cdots \phantom{.}
} \]
Further, $\mathrm{S}^{(a, b)}$ maps a morphism
\[ \xymatrix@C=30pt{ 
\: \phantom{\mathrm{S}^{(a, b)}} C: \ar@<10pt>[d]_{f} & \cdots \: \ar[r]  & C_{-1} \ar[r] \ar[d]^{f_{-1}} & C_0 \ar[r] \ar[d]^{f_{0}} & C_1 \ar[r] \ar[d]^{f_{1}} & \: \cdots \phantom{.} \\
\: \phantom{\mathrm{S}^{(a, b)}} D: & \cdots \: \ar[r]  & D_{-1} \ar[r] & D_0 \ar[r] & D_1 \ar[r] & \: \cdots \phantom{.}
} \]
to the following morphism of 2-complexes:
\[ \xymatrix@C=30pt{
\: \mathrm{S}^{(a, b)} C: \ar@<10pt>[d]_{\mathrm{S}^{(a, b)} f} & \cdots \: \ar[r]  & C_{-b} \ar[r] \ar[d]^{f_{-b}} & C_0 \ar[r] \ar[d]^{f_{0}} & C_a \ar[r] \ar[d]^{f_{a}} & \: \cdots \phantom{.} \\
\: \mathrm{S}^{(a, b)} D: & \cdots \: \ar[r]  & D_{-b} \ar[r] & D_0 \ar[r] & D_1 \ar[r] & \: \cdots .
} \]
As the classical homology
$$\HH^{(1, 1)}_0 = \ker d \big/ \im d : \Com{2}(\A) \to \A$$
is functorial, we derive that
\begin{displaymath}
\HH^{(a, b)}_0 = \HH^{(1, 1)}_0 \mathrm{S}^{(a, b)} 
\end{displaymath}
is functorial as well.

Denote by $\T^j$ the (invertible) translation functor that translates a complex $j$ steps ``to the left'',
that is position $i$ moves to position $i-j$. With this notation we have
\begin{equation}\label{eq713}
\HH^{(a, b)}_j = \HH^{(1, 1)}_0 \mathrm{S}^{(a, b)} \T^j.
\end{equation}
The claim of the proposition follows.
\end{proof}

The functors $\HH^{(a, b)}_j$ are called {\em Kapranov homology} functors.

\begin{proposition}\label{d-star}
The differential induces a natural transformation $\HH^{(a, b)}_{j} \to \HH^{(a, b)}_{j+1}.$
\end{proposition}

\begin{proof}
Apply $\HH^{(a, b)}_j$ to the morphism $\underline{d}_C : C \to \T^1 C$ given by 
\[ \xymatrix{ 
\cdots \: \ar[r]  & C_{-1} \ar[r]^d\ar[d]^{d} & C_0 \ar[r]^d \ar[d]^{d} & C_1 \ar[r]\ar[d]^d & \: \cdots \\
\cdots \: \ar[r]  & C_{0} \ar[r]^d & C_1 \ar[r]^d & C_2 \ar[r] & \: \cdots
} \]
Then, for each $C$, we have a sequence
\[ \xymatrix{
\cdots \: \ar[r] & \HH^{(a, b)}_{0}\T^{-1} C \ar[rr]^{\HH^{(a, b)}_{0} \underline{d}_{\T^{-1} C}} && \HH^{(a, b)}_{0}C \ar[rr]^{\HH^{(a, b)}_{0} \underline{d}_C} && \HH^{(a, b)}_{0}\T^1 C \ar[r] & \: \cdots
} \]
which is equal to the sequence
\begin{equation}\label{d-star-seq}
\xymatrix{
\cdots \: \ar[r] & \HH^{(a, b)}_{-1}C \ar[rr]^{\HH^{(a, b)}_{-1} \underline{d}_{C}} && \HH^{(a, b)}_{0}C \ar[rr]^{\HH^{(a, b)}_{0} \underline{d}_C} && \HH^{(a, b)}_{1}C \ar[r] & \: \cdots
}
\end{equation}
\end{proof}

The sequence \eqref{d-star-seq} is, in fact, an $M$-complex, where $M = \operatorname{min}(a, b)$.
Hence $\HH^{(a, b)}_{\_}$ defines a functor $$\HH^{(a, b)} : \Com{N}(\A) \to \Com{M}(\A).$$
In the classical setting of chain complexes, that is in the category $\Com{2}(\A) $, the functor
$\HH^{(1, 1)}$ maps $2$-complexes to $1$-complexes and, since $d = 0$ for $1$-complexes, 
the above property is not really visible.

\subsection{Total Kapranov homology}\label{homology.K}

As was shown in \cite{kapr91}, the Kapranov homology functors $\HH^{(a, b)}_j$ are connected by two families 
of commuting  natural transformations. We recall this construction here.

Let $C$ be an $N$-complex as in \eqref{eq2.1}.
The diagram \eqref{eq-inclusion} is a poset of inclusions. Therefore it commutes and  we may consider 
this diagram as the following morphism of sequences
\begin{equation}\label{eq-hom_i}
\xymatrix@C=23pt{ 
\cdots \: \ar[r]  & 0 \ar[r] \ar[d] & \im d^{N\text{-}1} \ar[r] \ar[d] & \im d^{N\text{-}2} \ar[r] \ar[d] 
& \cdots \ar[r] & \im d \ar[r] \ar[d] & 0 \ar[r] \ar[d] & \: \cdots \\
\cdots \: \ar[r]  & 0 \ar[r]        & \ker d      \ar[r]       & \ker d^2      \ar[r]        
& \cdots \ar[r] & \ker d^{N\text{-}1} \ar[r]        & 0 \ar[r]        & \: \cdots
}
\end{equation}
which we call $\alpha$. Since $\Seq(A)$ is abelian, $\alpha$ has a cokernel with objects being the homologies
\[ \xymatrix@C=21pt{ 
\cdots \: \ar[r]  & 0 \ar[r]        & \HH^{(1, N\text{-}1)}_{j}      \ar[r]^(.45){i_*}  
& \HH^{(1, N\text{-}2)}_{j}      \ar[r]^(.65){i_*}        
& ... \ar[r]^(.35){i_*} & \HH^{(N-1, 1)}_{j} \ar[r]        & 0 \ar[r]        & \: \cdots
} \]
which automatically defines the natural transformation $i_*$, the differential in the above cokernel sequence, between the corresponding homology functors.

The differential $d$ adds a new dimension to diagram \eqref{eq-hom_i} by inducing the vertical morphisms in the 
following commutative diagram:
\begin{equation}\label{eq-bi-im}
\xymatrix@!C=32pt{
0 \;\; \ar[r]  & \im d^{N\text{-}1} \ar[r] \ar[d]^d & \im d^{N\text{-}2} \ar[r] \ar[d]^d & ... \ar[r]  & \im d \ar[r] \ar[d]^d & 0 \ar[d]^0 & \\
& 0 \ar[r]        & \im d^{N\text{-}1} \ar[r] & ... \ar[r] & \im d^{2} \ar[r] & \im d \ar[r] & 0
}
\end{equation}
where the rows are coming from the top row of \eqref{eq-hom_i}. The same thing can be done with 
the bottom row of \eqref{eq-hom_i}:
\begin{equation}\label{eq-bi-ker}
\xymatrix@!C=32pt{
0 \;\; \ar[r]  & \ker d \ar[r] \ar[d]^d & \ker d^{2} \ar[r] \ar[d]^d & ... \ar[r] & \ker d^{N\text{-}1} \ar[r] \ar[d]^d & 0 \ar[d]^0 & \\
& 0 \ar[r]        & \ker d \ar[r] & ... \ar[r] & \ker d^{N\text{-}2} \ar[r] & \ker d^{N\text{-}1} \ar[r] & 0
}
\end{equation}
Indeed, diagrams \eqref{eq-bi-im} and \eqref{eq-bi-ker} both commute and are sequence of sequences
(we will call this bisequences in analogy with bicomplexes being complexes of complexes). 
Again, the morphism $\alpha$  in \eqref{eq-hom_i} defines a morphism of bisequences from 
the bisequence \eqref{eq-bi-im} to the bisequence \eqref{eq-bi-ker}. 
The cokernel  of this morphism of bisequences results in the following bisequence with components being 
Kapranov homology functors
\begin{equation}\label{eq-bi-hom}
\xymatrix@!C=32pt{
0 \;\; \ar[r]  & \HH^{(1, N\text{-}1)}_{j} \ar[r]^(.5){i_*} \ar[d]^{d^*} & \HH^{(2, N\text{-}2)}_{j} 
\ar[r]^(.62){i_*} \ar[d]^{d^*} & ... \ar[r]^(.38){i_*} & 
\HH^{(N\text{-}1, 1)}_{j} \ar[r] \ar[d]^{d^*} & 0 \ar[d]^0 & \\
& 0 \ar[r]        & \HH^{(1, N\text{-}1)}_{j+1} \ar[r]^(.62){i_*} & ... \ar[r]^(.38){i_*} & \HH^{(N\text{-}2, 2)}_{j+1} \ar[r]^(.5){i_*} & \HH^{(N\text{-}1, 1)}_{j+1} \ar[r] & 0
}
\end{equation}
where $d^*$ denotes the maps induced by $d$.

Define a new sequence with the objects
$$\mathbf{H}_n = \bigoplus_{2j + p = n} \HH^{(p, N-p)}_{j}$$
at position $n$
and the differential $\mathbf{H}_n \to \mathbf{H}_{n+1}$ given by adding up the $i_*$'s  and the $d^*$'s between the
corresponding summands. In \cite{kapr91} it is shown that this defines a functor 
$$\mathbf{H} : \Com{N}(\A) \to \Com{N-1}(\A),$$
the {\em total Kapranov homology functor}.

\subsection{Homology for Sequences}\label{s4}

Let $\A$ be an abelian category and $N$ a positive integer. Let $\T^1$ be the translation functor as defined above. 
Given an object $C$ in $\Seq(\A)$, the differential of $C$ gives rise to a morphism 
$$\underline{d}_C:C \longrightarrow \T^1 C$$
in $\Seq(\A)$ given by 
\begin{displaymath}
\xymatrix{
\dots\ar[r]^{d}&C_{-1}\ar[r]^{d}\ar[d]^{d}&C_{0}\ar[r]^{d}\ar[d]^{d}&C_{1}\ar[r]^{d}\ar[d]^{d}&\dots\\
\dots\ar[r]^{d}&C_{0}\ar[r]^{d}&C_{1}\ar[r]^{d}&C_{2}\ar[r]^{d}&\dots\\
}
\end{displaymath}
For simplicity we will omit  indices and, for example,  write $\underline{d}^2$ for $\underline{d}_{\T^1 C} \underline{d}_C$.

Let $\B$ be a full subcategory of $\Seq(\A)$ containing $\Com{N}(\A)$. 
Denote by  $\I$ the inclusion functor $\Com{N}(A) \hookrightarrow \B$.
For each object $C$ in $\B$  we have the diagram
\begin{equation}\label{eq2.2}
\xymatrix@C=30pt{ 
\ker \underline{d}^N_C \ar@<-1pt>@{^{(}->}[r]^(.6){k_C} & C \ar@<-1pt>[r]^(.45){\underline{d}^N} & \T^N C,
} 
\end{equation}
where $k_C$ denotes the inclusion.
Note that $\ker \underline{d}^N_C$ is an $N$-complex. As $\Com{N}(A) \hookrightarrow \B$, it follows that
the diagram \eqref{eq2.2} is, in fact, a diagram in $\B$. 

\begin{lemma}\label{lem201}
Let $C,D\in \Seq(\A)$ and $f:C\to D$ be a morphism. 
Then there is a unique $f':\ker \underline{d}^N_C\to  \underline{d}^N_D$ such that the following diagram commutes
\[ \xymatrix@C=30pt{ 
\ker \underline{d}^N_C \ar@{.>}[d]_{f'} \ar@<-1pt>@{^{(}->}[r]^(.6){k_C} & C \ar[d]^{f} 
\ar@<-1pt>[r]^{\underline{d}^N} & \T^NC \ar[d]^{\T^N(f)} \\
\ker \underline{d}^N_D \ar@<-1pt>@{^{(}->}[r]^(.6){k_D} & D \ar@<-1pt>[r]^{\underline{d}^N} & \T^ND
} \]
\end{lemma}

\begin{proof}
We have $\underline{d}^Nfk_C = \T^N(f)\underline{d}^Nk_C = 0$ since $\underline{d}^Nk_C = 0$. By the
universal property of kernels, there is a unique $f'$ such that $k_Df' = fk_C$. 
\end{proof}

Define $$\dd: \B \longrightarrow \Com{N}(A)$$ as follows: 
\begin{itemize}
\item a sequence $C$ is mapped to the $N$-complex $\ker \underline{d}^N_C$,
\item a morphism $f: C \to D$ is mapped  to $f'$ given by Lemma~\ref{lem201}.
\end{itemize}

\begin{lemma}\label{lem202}
We have that $\dd$ is a functor.
\end{lemma}

\begin{proof}
The composition $\dd(f)\dd(g)$ satisfies the universal property required for $\dd(fg)$ since the 
following diagram commutes:
\[ \xymatrix@C=30pt{ 
\ker \underline{d}^N_B \ar@{->}[d]_{\dd(g)} \ar@<-1pt>@{^{(}->}[r]^(.6){k_B} & B \ar[d]^{g} 
\ar@<-1pt>[r]^{\underline{d}^N} & \T^NB \ar[d]^{\T^N(g)} \\
\ker \underline{d}^N_C \ar@{->}[d]_{\dd(f)} \ar@<-1pt>@{^{(}->}[r]^(.6){k_C} & C \ar[d]^{f} 
\ar@<-1pt>[r]^{\underline{d}^N} & \T^NC \ar[d]^{\T^N(f)} \\
\ker \underline{d}^N_D \ar@<-1pt>@{^{(}->}[r]^(.6){k_D} & D \ar@<-1pt>[r]^{\underline{d}^N} & \T^ND
} \]
Because of the uniqueness, we deduce that $\dd(fg)$ and  $\dd(f)\dd(g)$ coincide. 
Similarly one shows that $\dd$ maps identities to identities.  
\end{proof}

\section{Derived Equivalence}\label{s5}

\subsection{An auxiliary adjunction}\label{s5.1}

\begin{proposition}\label{prop35}
The pair $(\I,\dd)$ is an adjoint pair of functors.
\end{proposition}

\begin{proof}
Directly from the definition we see that $\dd$ is the identity on $\Com{N}(\A)$, hence 
$$\dd \I = \Id_{\Com{N}(\A)}.$$
We claim that the unit of the adjunction is the identity natural transformation $$id: \Id_{\Com{N}(\A)} \longrightarrow \dd \I = \Id_{\Com{N}(\A)},$$ and  the counit of the adjunction is the inclusion 
$$k: \I \dd \longrightarrow \Id_{\B}$$
defined in Subsection~\ref{s4}.

Indeed, let $C,D\in\B$ and $f:C\to D$. Then $k_D \dd \I(f) = f k_C$ by construction. This means that 
$k$ is indeed a natural transformation. Further, we have $$k_{\I} \circ \I(id) = k_{\I} = id_{\I}$$ since the inclusion of a complex into itself is the identity. But also  $$\dd(k) \circ id_{\dd} = \dd(k) = id_{\dd}$$ 
because $k$ is monic and $k \circ \dd(k) = k \circ id_{\dd}$.
\end{proof}

\subsection{Localization of categories}\label{s5.2}

Recall that $\operatorname{CAT}$ denotes the metacategory of all categories, see \cite{mac}.
Given a category $\C$ and a subclass of morphisms $S \subset \Mor(\C)$, define the {\em localizing functor} $Q: \C \longrightarrow \C[S]$ as the universal functor in $\operatorname{CAT}$  that maps morphisms 
$f \in S$ to isomorphisms, that is, $Q: \C \longrightarrow \C[S]$ is a functor such that for any other functor 
$F : \C \to \D$ that maps morphisms in $S$ to isomorphisms, there is a unique, up to isomorphism,
functor $F' : \C[S] \to \D$  such that $F = F' Q$. The domain of $Q$, denoted $\C[S]$, is called the {\em localized category} for details.

It follows from the definition that if a localized category $\C[S]$ (with associated $Q$) exists, 
then it is unique up to isomorphism. The localized category can be constructed by adding 
formal inverses to morphisms in $S$, see \cite[III~\S{2}]{gelfand} for details.

\subsection{Main result}\label{s5.3}

\begin{theorem}\label{dereq}
Let $\C$ and $\D$ be two categories and $\G$ and $\F$ a pair of functors as follows:
\[ \xymatrix{ {\C} \ar@<2pt>[r]^{\G} & {\D} \ar@<2pt>[l]^{\F}. } \]
Assume we are given two natural transformations
$$ \eta :\Id_{\D} \longrightarrow \G\F , \qquad  \eps : \F\G \longrightarrow \Id_{\C} $$ and a pair of subclasses of morphisms
$$S \subset \Mor(\C) , \qquad T \subset \Mor(\D)$$
such that
\begin{displaymath}
\F(T) \subset S,\quad
\G(S) \subset T,\quad 
\{\eps_{X}\}_{X \in \C} \subset S,\quad
\{\eta_{X}\}_{X \in \D} \subset T.
\end{displaymath}
Then the localized categories $\C[S]$ and $\D[T]$ are equivalent. Moreover, the same holds if we arbitrarily change the direction of $\eta$ or $\eps$.
\end{theorem}

\begin{proof}
Let $\Q : \C \to \C[S]$ and $\R : \D \to \D[T]$ be the corresponding localizing functors. Since $\R\G$
maps morphisms from $S$ to isomorphisms, there is a unique $\G'$ such that $\R\G = \G'\Q$.
Similarly, there is a unique $\F'$  such that $\Q\F = \F'\R$, that is the following diagrams commute 
\[ \xymatrix{
\C \ar[d]_{\Q} \ar[r]^{\G} & \D \ar[d]^{\R} &
\C \ar[d]_{\Q} & \D \ar[l]_{\F} \ar[d]^{\R} \\
\C[S] \ar[r]^{\G'} & \D[T] &
\C[S] & \D[T] \ar[l]_{\F'}
} \] \\
We have a natural transformation $\Q(\eps): \Q\F\G =_{\HH} \F'\G'\Q \to \Q$.
But $\Q$ acts as the identity on objects, so for each object $X$ of $\C[S]$ we have an isomorphism
$$\eps_X: \F' \G' X \longrightarrow X$$
in $\C[S]$ because $\eps_X$ lies in $S$.
To show that $\eps$ defines a natural transformation $\F'\G' \to \Id_{\C[S]}$, it is enough to check that $\eps_X$ commutes with generators in $\C[S]$.

We already know that $\eps_X$ commutes with the morphisms of $\C$. Since morphisms from $S$
are invertible in $\C[S]$, it follows that $\eps_X$ commutes with the inverses of all morphisms in $S$
(indeed, $xs=ty$ implies $t^{-1}x=ys^{-1}$ if $s$ and $t$ are invertible).
Since all $\eps_X$ are invertible in $\C[S]$ and $\eps$ is a natural transformation, similarly the inverses $\eps^{-1}_X$
consitute a natural transformation, call it $\eps^{-1}$, and this natural transformation is the inverse of
$\eps$. Then $$\eps \eps^{-1} = id_{\Id_{\C}} \mbox{ and } \eps^{-1} \eps = id_{\F'\G'}$$ which proves that the
functors $\Id_{\C}$ and $\F'\G'$ are isomorphic. Similarly one shows that 
$\Id_{\D}$ and $\G'\F'$ are isomorphic. This proves the first claim of the theorem.

Since the direction of $\eps$ (resp. $\eta$) did not matter in the argument above, the same argument can be applied
if we arbitrarily change the direction of $\eps$ or $\eta$.
\end{proof}

\subsection{Derived categories}\label{s5.4}

Let $\C$ and $\D$ be two categories and $\HH : \C \to \D$ be a functor. A morphism $f$ of $\C$ is
called an \textit{$\HH$-quasi-isomorphism} if $\HH(f)$ is an isomorphism.

Further, we define the corresponding {\em derived category} $\DD_{\HH}(\C)$
as the category $\C$ localized with respect to all $\HH$-quasi-isomorphisms, that is
\begin{equation}\label{quasi-iso}
\DD_{\HH}(\C) := \C[S] \,\,\text{ where }\,\, S = \{f \in \Mor(\C) : \HH(f) \text{ is an isomorphism}\}.
\end{equation}
For example, for an abelian category $\A$ the usual derived category $\DD(\A)$ is defined as
$\DD_{\HH^{(1,1)}}(\Com{2}(\A))$, see \cite[III~\S{2}]{gelfand}. Similarly one can define derived categories
of sequences or $N$-complexes choosing one's favorite homology functor.

\subsection{Derived equivalence}\label{s5.5}

\begin{theorem}\label{theorem_eq}
Let $\A$ be an abelian category, $N\in\{1,2,3,\dots\}$ and $\B$ be a full subcategory of $\Seq(\A)$ containing
$\Com{N}(\A)$. Then, given any category $\C$ and any functor $\E : \Com{N}(\A) \to \C$, there is an 
equivalence of categories
$$\DD_{\HH}(\B) \simeq \DD_{\HH}(\Com{N}(\A))$$
where $\HH = \E \dd$.
\end{theorem}

\begin{proof}
Recall that $\dd$ acts as the identity functor on $\Com{N}(\A)$.
We want to apply Theorem~\ref{dereq}.
The adjoint pair $(\I, \dd)$
\[ \xymatrix@C=50pt{ {\B} \ar@<2pt>[r]^(.4){\dd} & {\Com{N}(\A)} \ar@<2pt>[l]^(.6){\I}, } \]
given by Propositiuon~\ref{prop35}, and the corresponding adjunction morphisms 
provide  the functors and, respectively, the natural transformations as mentioned in the formulation of 
Theorem~\ref{dereq}.
The subclasses of morphisms to be localized are the $\HH$-quasi-isomorphisms as described in \eqref{quasi-iso}.
The functor $\dd : \B \to \Com{N}(\A)$ maps quasi-isomorphisms to quasi-isomorphisms because $\HH \dd = \HH$.
The functor $\I : \Com{N}(\A) \to \B$ maps quasi-isomorphisms to quasi-isomorphisms because 
$\HH \I = \E \dd \I = \E = \HH$.

Finally, let us check that the natural transformations $id$ and $k$ consist of quasi-isomorphisms. 
For $id$ it is obvious and for $k$ it holds because $\dd(k) = id$. Therefore we may apply Theorem~\ref{dereq}
which yields the desired statement.
\end{proof}

The following proposition deals with the case $\E = \HH^{(a, b)}$ in Theorem \ref{theorem_eq} which is of particular interest.

\begin{proposition}\label{homology_reformulation}
Let $N = a + b$. Then the functor $\HH^{(a, b)} \dd$ can be written in the form
\begin{equation*}
\ker d^a \big/ \left( \ker d^a \cap \im d^b \right) \text{ or } \left( \ker d^a + \im d^b \right) \big/ \im d^b.
\end{equation*}
\end{proposition}
\begin{proof}
The two expressions are equal by the second isomorphism theorem.
We show this for the left hand side expression.
Let $C$ be a sequence and denote the differential of $C$ and $\dd C$ by $d$ and $\delta$, respectively. Then $\ker d^a = \ker \delta^a$ because $\ker d^a \subset \ker d^N$.
In general $\im d^b \neq \im \delta^b$ but we have
$$\ker d^a \cap \im d^b = \ker \delta^a \cap \im \delta^b = \im \delta^b$$
because $a + b = N$.
\end{proof}

\begin{corollary}
Let homology on $\Com{N}(\A)$ be given by $\HH = \ker d / \left( \ker d \cap \im d \right)$.
Then there is an equivalence of categories
$$\DD_{\HH}(\Com{N}(\A)) \simeq \DD_{\HH}(\Com{2}(\A)).$$
\end{corollary}

\renewcommand{\dd}{[\ker d^2]}
\begin{proof}
By Proposition \ref{homology_reformulation} we have $\HH = \HH^{(1, 1)} \dd$.
The statement follows from Theorem \ref{theorem_eq} with $\B = \Com{N}(\A)$.
\end{proof}
\renewcommand{\dd}{[\ker d^N]}

\begin{corollary}
There is an equivalence of categories
$$\DD_{\dd}(\B) \simeq \Com{N}(\A).$$
\end{corollary}

\begin{proof}
This follows from  Theorem \ref{theorem_eq} by taking $\C=\Com{N}(\A)$ and
$\E = \mathrm{Id}_{\Com{N}(\A)}$ as under this choice we have  $\DD_{\HH}(\Com{N}(\A)) = \Com{N}(\A)$
since $\Com{N}(\A)$ has the universal property required for the localized category.
\end{proof}

\section{Projective Resolutions}\label{s6}

We continue to work with categories $\A$, $\Seq(\A)$ and $\Com{N}(\A)$ as described above.
Denote by $\R_N : \Com{2}(\A) \to \Com{N}(\A)$ the functor mapping a chain complex 
\[ C:\quad\xymatrix{ 
\cdots \: \ar[r]^d  & C_{-1} \ar[r]^d & C_0 \ar[r]^d & C_1 \ar[r]^d & C_2 \ar[r]^d & \: \cdots
} \]
to the $N$-complex $\R_N C$ where
$$(\R_N C)_{Nj} := C_{2j}\quad\text{ for }\quad j \in \Z,$$
and the odd degrees of $C$ are repeated $N - 1$ times as follows:
\[ \xymatrix@=20pt@R=1pt{ 
\R_N C:&\cdots \: \ar[r]^1  & C_{-1} \ar[r]^d & C_0 \ar[r]^d & C_1 \ar[r]^1 & \cdots \ar[r]^1 & C_1 \ar[r]^d & C_2 \ar[r]^d & \: \cdots \\
\text{degree}&&{-1}&0&1&&N-1&N&
} \]
The action of  $\R_N$ on morphisms is defined accordingly.

\begin{proposition}\label{preserve_H-quasi-iso}
For $a,b\in\{1,2,3,\dots\}$ the functor $\R_{(a+b)}$ maps ordinary quasi isomorphisms of chain complexes to 
$\HH^{(a, b)}$-quasi isomorphisms of $(a+b)$-complexes.
\end{proposition}

\begin{proof}
This follows from the definitions and \eqref{eq713} by a direct computation.
\end{proof}

Recall that an abelian category $\A$ is said to have {\em enough projectives} provided that each object in $\A$
is a quotient of a projective object. A natural analogue of this standrd notions for sequences over $\A$ looks
as follows:

\begin{definition}
Let $\B$ be  a full subcategory of $\Seq(\A)$, $\C$ any category and $\HH : \B \to \C$ a functor.
An object $P \in \B$ is said to be a \textit{$\HH$-projective resolution} of the object $X \in \A$ if $P_i$ is equal to zero for $i \geq 1$, $P_i$ is projective for $i \leq 0$ and there is an $\HH$-quasi isomorphisms in $\B$ as follows
\[ \xymatrix@C=30pt{
\cdots \: \ar[r]  & P_{-2} \ar[r] \ar[d] & P_{-1} \ar[r] \ar[d] & P_0 \ar[r] \ar[d] & 0 \ar[r] \ar[d] & 0 \ar[r] \ar[d] & \: \cdots \phantom{.}\\
\cdots \: \ar[r]  & 0 \ar[r] & 0 \ar[r] & X \ar[r] & 0 \ar[r] & 0 \ar[r] & \: \cdots .
} \]
\end{definition}

\begin{corollary}\label{enough_proj}
Let $a,b\in\{1,2,3,\dots\}$. If $\A$ has enough projectives, then every object in $\A$, considered as an object in 
$\Com{a+b}(\A)$ in the obvious way, admits an $\HH^{(a, b)}$-projective resolution.
\end{corollary}

\begin{proof}
If some object $P\in\Com{2}(\A)$ is an $\HH^{(1, 1)}$-projective resolution of $X \in \A$
(i.e. $P$ is a projective resolution of $X$ in the usual sense), then, by Proposition~\ref{preserve_H-quasi-iso},
the object $\R_{(a+b)} P\in\Com{a+b}(\A)$ 
is an $\HH^{(a, b)}$-projective resolution of $X$.
\end{proof}

The following statement gives a kind of a ``lower bound on $N$'' for existence of projective resolutions.

\begin{theorem}\label{proj}
Let $X \in \A$ be a non-projective object and $P \in \Com{a+b}(\A)$ an $\HH^{(a, b)}$-projective resolution of $X$.
Then $P \notin \Com{a+b-1}(\A)$.
\end{theorem}

\begin{proof}
We need to find a position at which the differential $d$ of $P$ satisfies
$$d^{a + b - 1} \neq 0.$$

For the case $a = 1$ we have $$\HH^{(a, b)}_0(P) = P_0 / d^{b} P_{-b} \simeq X$$ but $d^{b} P_{-b} = d^{a+b-1} P_{-b} \neq 0$ since $X$ was not projective.

Now assume $a \geq 2$. Homology of $P$ at position $1-a$ is
$$\HH^{(a, b)}_{1-a}(P) = P_{1-a} / d^{b} P_{1-a-b} = 0$$
and hence
$$d^{b} : P_{1-a-b} \to P_{1-a}$$
is epi. If $1 - a < -b$, homology of $P$ at position $1-a+b$ is again zero and hence
$$d^{b} : P_{1-a} \to P_{1-a+b}$$
is epi.
We may now repeat this $n$ more steps until $1-a+nb \geq -b$. But, similarly to the case of $a = 1$, 
$d^{b} : P_{-b} \to P_0$ is non-zero.
This means, in particular, that $d^{-1+a-nb} : P_{1-a+nb} \to P_0$ is non-zero and 
$$d^{-1+a-nb} (d^b)^{n+1} = d^{a + b - 1} : P_{a + b - 1} \to P_0$$
is non-zero. The claim follow.
\end{proof}

The combination of Theorems~\ref{theorem_eq} and \ref{proj} suggests that the category $\Com{N}(\A)$ 
is  ``homologically optimal'' with respect to (generalized) Kapranov homology functors in the sense that, 
on the one hand,  it is big enough so that the derived category inherits the property of having projective 
resolutions and, on the other hand, it is small enough in the sense the derived category of any bigger 
category produces an equivalent category.
\vspace{1cm}

{\bf Acknowledgements.}
I thank Volodymyr Mazorchuk for his advising.

\end{document}